\documentclass[a4paper,11pt]{article}
\usepackage[english]{babel}
\usepackage{theorem}
\usepackage{amssymb}
\usepackage{amsfonts}
\usepackage{amsmath}
\usepackage[all]{xy}
\usepackage{geometry}
\usepackage{ae,aecompl}
\usepackage{eurosym}

\newtheorem{thm}{Theorem}[section]
\newtheorem{prop}[thm]{Proposition}
\newtheorem{lem}[thm]{Lemma}
\newtheorem{cor}[thm]{Corollary}

\makeatletter
\begingroup
\gdef\th@upshape{\normalfont
  \def\@begintheorem##1##2{%
        \item[\hskip\labelsep \theorem@headerfont ##1\ ##2.]}%
\def\@opargbegintheorem##1##2##3{%
   \item[\hskip\labelsep \theorem@headerfont ##1\ ##2\ (##3).]}}
\endgroup
\makeatother

\theoremstyle{upshape}

\newtheorem{defn}[thm]{Definition}

\title{{\bf The Noncommutative K\"ahler Geometry of the Standard Podle\'s Sphere}}
\author{R\'{e}amonn \'{O} Buachalla}
\def\bal#1\eal{\begin{align}#1\end{align}}
\def\bas#1\eas{\begin{align*}#1\end{align*}}
\def\bit{\begin{itemize}}
\def\eit{\end{itemize}}
\def\ed{\end{document}}

\def\e{\varepsilon}

\def\k{\kappa}

\def\t{\tau}
\def\w{\omega}

\def\Om{\Omega}

\def\del{\partial}
\def\adel{\ol{\partial}}
\def\DEL{\Delta}

\def\bC{{\mathbf C}}
\def\bN{{\mathbf N}}

\def\bZ{{\mathbf Z}}

\def\E{{\cal E}}

\def\exd{\mathrm{d}}
\def\demo{\noindent \emph{\textbf{Proof}.\ }~}

\def\frame{\textrm{frame}_M}

\def\id{\mathrm{id}}
\def\ker{\mathrm{ker}}

\def\spn{\mathrm{span}}

\def\hol{^{(1,0)}}
\def\ahol{^{(0,1)}}

\def\inv{^{-1}}

\def\oby{\otimes}

\def\wed{\wedge}
\def\sseq{\subseteq}

\def\wt{\widetilde}

\def\ol{\overline}
\def\la{\left\langle}
\def\\la{\left\langle}
\def\ra{\right\rangle}
\def\>{\right\rangle}

\def\mto{\mapsto}

\def\qed{\hfill\ensuremath{\square}\par}
\def\qf3{\bC_q[F_3]}

\def\csu2{\bC_q[SU_2]}
\def\cs2{\bC_q[S^2]}

\def\cp1{\bC_q[\bC P^{1}]}
\def\ccp1{\bC P^{1}}

\def\usl2{\mathcal{U}(\mathfrak{sl}(2))}

\def\ws2{\Om^1_q(S^2)}

\def\m1{_{(-1)}}
\def\0{_{(0)}}
\def\1{_{(1)}}
\def\2{_{(2)}}
\def\3{_{(3)}}
\def\4{_{(4)}}
\def\5{_{(5)}}
\def\hol{^{(1,0)}}
\def\ahol{^{(0,1)}}
\def\asth{\ast_H}
\def\lapd{\DEL_{\exd}}
\def\lapdel{\DEL_{\del}}
\def\lapadel{\DEL_{\adel}}

\def\cvect{{}_{\bC} \hspace{-.030cm}\mathcal{M}}

\def\mm{{}_M \hspace{-.030cm}\mathcal{M}}

\def\alg{algebra~}
\def\algn{algebra}

\def\nccg{noncommutative complex geometry~}
\def\nccgn{noncommutative complex geometry}

\def\nc{noncommutative~}

\def\st{such that~}

\def\wrt{with respect to~}

\def\uqsl2{U_q(\frak{sl_2)}}
\def\tuqsl2{\wt{U}_q(\frak{sl}_2)}
\def\tu1sl2{\wt{U}_1(\frak{sl}_2)}

\begin{document}

\maketitle

\begin{abstract}
Building on the now established presentation of the standard Podle\'s sphere as an example of a noncommutative complex structure, we investigate how its classical  K\"ahler geometry behaves under $q$-deformation. Discussed are noncommutative versions of Hodge decomposition, Lefschetz decomposition, the K\"ahler identities,  and the refinement  of  de Rham cohomology by Dolbeault cohomology.%, and the Schr\"odinger--Lichnerowicz formula.
\end{abstract}

%%%%%%%%%%%%%%%%%%%%%%%%%%%%%%%%%%%%%%%%%%%%%%%%%%%%%%%%%%%%%%%%%%%%%%%%%%%%%%%%%%%%%%%%
\section{Introduction}
%%%%%%%%%%%%%%%%%%%%%%%%%%%%%%%%%%%%%%%%%%%%%%%%%%%%%%%%%%%%%%%%%%%%%%%%%%%%%%%%%%%%%%%

%%%%%%%  (A) Podles and Non-Comm Geom

Since its introduction, the Podle\'s sphere  \cite{PodlesSphere}  has served as an example of central importance for many areas of noncommutative geometry. We cite the role it has played in the theory of covariant differential calculi \cite{Podcalc, PodcalcClass, Herm}, quantum principal bundles \cite{qmonop}, quantum frame bundles \cite{Maj, Maj1}, cyclic cohomology \cite{MNW, Had}, and spectral triples \cite{DS, RS}.

%.... these have all been subsequently generalised . . . HK, Meyer, MMF1, uli, HK, Dirac CPn, MMF

\bigskip

%%%%%%%  (B) Podles and Non-Comm C-Geom
In recent years, the Podle\'s sphere has assumed a similarly important role in the newly emerging field of \nccgn. %Given the fundamental role of the two sphere in classical complex geometry, this should not be so surprising. 
The existence for the Podle\'s sphere of a $q$-deformed Dolbeault double complex was discovered independently by Majid, and by Heckenberger and Kolb.  In \cite{Maj} it was arrived at using a frame bundle approach,  while in \cite {HK} it emerged from a classification of the covariant first order differential calculi of the irreducible quantum flag manifolds. A definition of \nc complex structure would later be introduced in \cite{KLVS} in order to formalize the properties of this $q$-Dolbeault complex. A subsequent more comprehensive version of this definition would appear in \cite{EBPS}, following which a third version, tailored for quantum homogeneous spaces, was introduced by the author in \cite{MMF2}. Aspects of the noncommutative Hermitian geometry of the Podle\'s sphere have also been investigated. In \cite{LZ1} and \cite{Maj} there appeared quantum versions of the $2$-sphere's covariant Hermitian metric, along with an associated Hodge operator.

\bigskip

%%%%%%%  (C) Podles and Non-Comm K-Geom

Classically, $S^2_q$ is not just a Hermitian manifold, but a K\"ahler manifold. The goal of this paper is to build upon the \nc complex and Hermitian constructions outlined above, and to propose the Podle\'s sphere as a prototypical example of a \nc K\"ahler structure. To justify this proposal we establish noncommutative versions of Hodge decomposition, Lefschetz decomposition, the K\"ahler identities,  and the refinement of  de Rham cohomology by Dolbeault cohomology. An abstract definition for \nc K\"ahler structure will appear in \cite{MMF4}, along with general proofs of the results established in this paper, and applications thereof to the quantum projective spaces. %We should also note that elements of  K\"abler structure for $\cs2$ have already appeared in \cite{oldkahler}, and it would be interesting to link them with the work of this paper

\bigskip

%%%%%%%%%%%%%%%% E Achoimre

%%%  2
The paper is organised as follows: In section 2 we will recall the definition of the Podle\'s sphere as a quantum homogeneous space of $\bC_q[SU_2]$; its quantum Dolbeault double complex; as well as its Hermitian metric and associated Hodge operator.

%%%  3
In section 3, we will show that the Hodge decompositions of $\Om^1(S^2)$ \wrt $\exd,\del$, and $\ol{\del}$ carry over directly to the quantum setting. This allows us to show that the dimensions of the  cohomology groups of $\Om^{\bullet}_q(S^2)$ have classical lower bounds. %%% Has the cohomology been calculated before

%%%  4
In section 4, we will introduce natural quantum analogues of the Lefschetz and dual Lefschetz operators. Moreover, we will show that these induce a direct generalisation of the Lefschetz decomposition.

%%%  5
In Section 5, we will show that the classical K\"ahler identities for the two sphere carry over to the quantum setting undeformed. Using this result we can then easily conclude that Dolbeault cohomology is a refinement of de Rham cohomology, generalising another important result of classical K\"ahler theory.

%%% 6
%Finally, in Section 6, we conclude from direct calculation that a version of the classical Schr\"odinger--Lichnerowicz formula is satisfied for $\bC_q[S^2]$. This result is expected to be very useful in . . .  \nc Sobolev theory.

\section{Preliminaries}

In this section we fix notation and recall the definitions, constructions, and results from the basic theory of the Podle\'s sphere. References are provided where proofs or basic details are omitted. 

\subsection{The Podle\'s Sphere}

In this subsection we will first recall some basic facts about general theory of faithfully flat quantum homogeneous spaces. We will then consider the standard Podle\'s sphere as an example.

\subsubsection{Faithfully Flat Quantum Homogeneous Spaces}

Let $H,G$ be two Hopf $*$-\algn s, and let us denote the coproducts, counits, antipodes, and $*$-maps of both by $\DEL$, $\e$,$S$, and $*$ respectively. Moreover, for $\pi: G \to H$ a Hopf  $*$-\alg map,   let us write $\DEL_{\pi} := (\id \oby \pi) \circ \DEL$. We denote the coinvariant $*$-sub\alg of $\DEL_{\pi}$ by $G^H$, that is
\bas 
G^H := \{g \in G \,|\, \DEL_{\pi}(g) = g \oby 1\}.
\eas
We call an \alg of the form $G^H$ a {\em quantum homogeneous space}.  An important fact is that every quantum homogeneous space has a canonical left $G$-coaction  $\DEL_L: M \to G \oby M$, induced in the obvious way by the coproduct of $G$.

We say that $G$ is a {\em faithfully flat} module over $M$ if the {\em tensor product functor}  $G \oby_M -:\mm \to \cvect$, from the category of left $M$-modules to the category of complex vector spaces, preserves and reflects exact sequences. 

Let us now explain why faithful flatness is important to us: For a quantum homogeneous $M = G^H$, let $\E$ be an $M$-bimodule endowed with a left $G$-coaction $\DEL_L$, satisfying the compatibility condition
\bas
\DEL_L(m e m') =  m\1 e\m1 m'\1 \oby m\2 e\0 m'\2, & &( \text{for all ~} m,m' \in M, e \in \E).
\eas
Moreover, let us denote by $\Phi(\E)$, the right $H$-comodule $\E/(M^+\E)$ with coaction $\DEL_R(\ol{e}) =\ol{e\1} \oby \pi(S(e\2))$. Now if $G$ is a faithfully flat module over $M$, It follows from a result of Takeuchi \cite{Tak}, that we have an isomorphism
\bas
\frame: \E \to (G \oby \Phi(\E))^H, & & e \mto e\m1 \oby \ol{e\0}
\eas
For a more in depth description of this important result, see \cite{MMF2} and references therein.

\subsubsection{The Standard Podle\'s Sphere}

In this paper, the quantum homogeneous space we will be working with is the Podle\'s sphere. Let us begin recalling its definition by recalling the well-known choice of $G$ in this case. For $q \in \bC^\times$, the $\bC$-\alg $\bC_q[SL_2]$ is generated by the four elements $a,b,c,d$, subject to the relations 
\bas
ab = qba, & & ac = qca, & & bc = cb,\\
ad - da - qbc, & & ad - qbc - 1.
\eas
It can be given a coalgebra structure with a coproduct uniquely  determined by
\bas
\DEL(a) = a \oby a + b \oby c, & & \DEL(b) = a \oby b + b \oby d, \\
 \DEL(c) = c \oby a + d \oby c, & & \DEL(d) = c \oby b +  d \oby d,
\eas
and a counit uniquely determined by $\e(a) = \e(d) = 1$, and $\e(b) = \e(c) = 0$. This coalgebra structure is easily seen to be extendable  to a Hopf \alg structure. The corresponding antipode $S$ satisfies the relations 
$S(a) = d, S(b) = -q\inv c, S(c) = -q b,$ and  $S(d) = a$. Finally,  $\bC_q[SL_2]$ can be given a Hopf $*$-algebra structure uniquely determined by $\ast(a) = d, \ast(b) = -q b, \ast(c) = -q \inv c,$ and $\ast(d) = a$. When  $\bC_q[SL_2]$ is endowed with this $\ast$-structure, we denote it by $\bC_q[SU_2]$. Note that, for any $f \in \bC_q[SU(2)]$, we will use $\ast(f)$ and $f^*$ interchangeably.

Turning now to our choice for $H$, we recall that  $\bC[U_1]$  is  the commutative \alg over $\bC$ generated by $t$ and $t\inv$, subject to the obvious relation $t t \inv = t\inv t = 1$. It has a Hopf $*$-algebra structure uniquely determined by  $\DEL(t) = t \oby t$; $\e(t) = 1$; $S(t) = t\inv$; and $\ast(t) = t\inv$. 

Finally, we come to the question of a map from $G$ to $H$, and choose $\pi: \bC_q[SU(2)] \to \bC[U_1]$ to be the unique Hopf algebra map determined by
\bas
\pi(a) = t\inv, & & \pi(d) = t, & & \pi(b) = \pi(c) = 0.
\eas
We call the corresponding coinvariant sub\alg $\bC_q[SU_2]^{\bC[U_1]}$  the {\em Podle\'s sphere}, and denote it  by $\bC_q[S^2]$. As a unital algebra it is generated by the elements $b_- := ab, b_0:= bc, b_+ := cd$. Moreover, $\csu2$ is a faithfully flat module over $\cs2$ \cite{MulSch}.

We should note that the coaction $\DEL_\pi$ induces a $\bZ$-grading on $\bC_q[SU_2]$, where, for example, we have
\bas
\deg(a) = \deg(c) = -1, & & \deg(c) = \deg(d) = 1.
\eas
Clearly, $\bC_q[S^2]$ is the degree-$0$ part of  $\bC_q[SU_2]$. More generally, we denote the degree-$k$ part by $\E_k$.

\subsection{The Noncommutative Complex Geometry of the Podle\'s Sphere}

We will now recall the \nccg of the Podle\'s sphere referred to in the introduction. We begin by recalling the basic definitions of \nccg in general, and then move on to the specific case of the Podle\'s sphere.

\subsubsection{Noncommutative Complex Structures}

A pair $(\Om^\bullet,d)$ is called a {\em differential \algn} if $\Om^\bullet = \bigoplus_{k \in \bN_0} \Om^k$ is an $\bN_0$-graded \algn, and $d$ is a degree $1$ map \st  $d^2=0$, and for which  the {\em graded Liebniz rule} is satisfied
\bas
\exd(\w \wed \nu) = \exd \w \wed \nu + (-1)^k \w \wed \exd \nu, & & (\w \in \Om^k, \nu \in \Om^\bullet).
\eas
 A {\em total differential calculus} over an \alg $A$ is a differential 
\alg $(\Om(A),\exd)$, \st $\Om^0=A$, and
%\begin{align*} %\label{pofu}
$\Om^k = \spn_{\bC}\{a_0\exd a_1\wed  \cdots \wed \exd a_k \,|\, a_0, \ldots, a_k \in A\}$. 
We call a differential calculus $(\Om^\bullet, \exd)$ over a $*$-\alg $A$ a {\em total $*$-differential calculus}, if the involution of $A$ extends to an involutive conjugate-linear map $*$ on $\Om^\bullet$, for which $(\exd \w)^* = \exd \w^*$, for all $\w \in \Om$, and
\[
(\w_p\w_q)^*=(-1)^{pq}\w_q^*\w_p^*, \qquad \text{ (for all }\w_p \in
\Om^p,~ \w_q \in \Om^q).
\]

\begin{defn}
An {\em almost complex structure} for a total $*$-differential calculus  $\Om^{\bullet}(A)$ over a $*$-\alg $A$, is an $\bN^2_0$-\alg grading $\bigoplus_{(p,q)\in \bN^2_0} \Om^{(p,q)}$ for $\Om^{\bullet}(A)$ such that, for all $(p,q) \in \bN^2_0$: 
\begin{enumerate}
\item \label{compt-grading}  $\Om^k(A) = \bigoplus_{p+q = k} \Om^{(p,q)}$;
\item  the wedge map restricts to isomorphisms
\bal \label{wedge-cond} 
\wed:\Om^{(p,0)} \oby_A \Om^{(0,q)} \to \Om^{(p,q)}, & & \wed: \Om^{(0,q)} \oby_A \Om^{(p,0)} \to \Om^{(p,q)};
\eal
\item  \label{star-cond} $*(\Om^{(p,q)}) = \Om^{(q,p)}$.
\end{enumerate}
We call an element of $\Om^{(p,q)}$ a {\em  $(p,q)$-form}.
\end{defn}

We say that a total differential calculus $\Om^\bullet(M)$  over a quantum homogeneous space $M = G^H$ is {\em covariant} if the coaction $\DEL_L$ extends to a left coaction on $\Om^\bullet(M)$ \st 
\bas
\DEL_L \circ \exd = (\id \oby \exd) \circ \DEL_L.
\eas
Moreover, we say that an almost complex structure $\Om^\bullet(M)$ is {\em left-covariant} if we have 
\bas
\DEL_L(\Om^{(p,q)}) \sseq  G \oby \Om^{(p,q)}, & & (\text{for all } (p,q) \in \bN^2).
\eas

Directly generalising the classical definition, we say that an almost-complex structure $\Om^{(\bullet,\bullet)}$ is {\em integrable} if $\exd(\Om^{(1,0)}) \sseq \Om^{(2,0)} \oplus \Om^{(1,1)}$, or equivalently, if $\exd(\Om^{(0,1)}) \sseq \Om^{(1,1)} \oplus \Om^{(0,2)}$. The assumption of integrability has some very useful consequences.

\begin{lem}
If an almost complex structure $\bigoplus_{(p,q)\in \bN_0^2} \Om^{(p,q)}$ is integrable, then
\begin{enumerate}
\item $\exd = \del + \ol{\del}$;
\item $(\bigoplus_{(p,q)\in \bN^2}\Om^{(p,q)}, \del,\ol{\del})$ is a double complex;
\item $\del(a^*) = (\ol{\del} a)^*$, and $\ol{\del}(a^*) = (\del a)^*$,   for all $a \in A$;
\item both $\del$ and $\ol{\del}$ satisfy the graded Leibniz rule.
\end{enumerate}
\end{lem}

Finally, we come to the definition of cohomology groups:  Just as in the classical case, the {\em de Rham cohomology groups} of $\Om^\bullet_q(S^2)$ are 
$
H^{k}(S^2) = \ker(\exd|_{\Om^{k}})/(\exd|_{\Om^{k-1}}).
$
Moreover, the {\em holomorphic} and {\em anti-holomorphic Dolbeault cohomology groups} are respectively defined by
\bas
H^{(a,b)}_{\del}(S^2) = \ker(\del|_{\Om^{a,b}})/\text{im}(\del|_{\Om^{a-1,b}}), & & H^{(a,b)}_{\ol{\del}}(S^2) = \ker(\ol{\del}|_{\Om^{a,b}})/\text{im}(\ol{\del}|_{\Om^{a,b-1}}).
\eas

\subsubsection{A Noncommutative Complex Structure for the Podles Sphere}

%%%%%%%%%%%%%%%%%%%%%%%%
In the case of the  Podle\'s sphere, there exists only one covariant total differential calculus  $\Om^\bullet_q(S^2)$ whose $1$-forms have classical dimension and which admits a covariant almost complex structure. Moreover,  $\Om^\bullet_q(S^2)$ admits only one such almost complex structure,  and it is integrable.

While we will not give a complete description of the calculus here, we will need to recall some important facts. Firstly, let us denote 
\bas
V\hol := \Phi(\Om\hol), & & V\ahol := \Phi(\Om\ahol).
\eas
Both $V\hol$ and $V\ahol$ are $1$-dimensional. We choose a basis for both according to 
\bas
e^+ := \ol{\exd b_+}, & & e^- := \ol{\exd b_-}.
\eas
The right $\bC[U_1]$-comodule structures of each are determined by  $\DEL_R(e^\pm) = e^\pm \oby t^{\pm 2}$. This immediately implies that 
\bas
\Om\hol = \E_2 \oby e^+, & & \Om\ahol = \E_{-2} \oby e^-.
\eas 
Now $\Phi(\Om^2_q(S^2))$ is a $1$-dimensional vector space, for which we choose $\t := \ol{\exd b_+ \wed \exd b_-}$ as a basis. It is easily shown that $\DEL_R(\t) = \t \oby 1$, and so we have that   
\bas
\Om^2_q(S^2) \simeq \bC_q[S^2] \oby \t.
\eas
Finally, for all $k > 2$, we have that  $\Om^k_q(S^2) = 0$.
For $f \in \E_{}$, $g \in \E_{}$, the multiplication $\wed$ in $\Om^\bullet_q(S^2)$ is uniquely determined by $(f \oby e^\pm) \wed (g \oby e^\pm) = 0$, and 
\bas
(f \oby e^+) \wed (f \oby e^-) = fg \t, & & (g \oby e^-) \wed (g \oby e^+) = -q^{2} gf \t.
\eas

\subsection{The Noncommutative Hermitian Geometry of the Podle\'s Sphere}

Let us start this subsection by constructing a sesqui-linear  map
\bas
g:\Om^1_q(S^2) \oby_{\bC_q[S^2]} \Om^1_q(S^2) \to \bC_q[S^2],
\eas
which we call the {\em metric} of $\Om^1_q(S^2)$: For $f \in \E_{2}$, and $h \in \E_{-2}$, we define $g$  to be the unique sesqui-linear mapping for which
\bas
g(f e^+ \oby h e^-) = fh^* , & & g(he^- \oby f e^+) =  q^{\pm 2}hf^*,  & & g(f e^{+} \oby f e^{+}) =  g(h e^{-}\oby h e^-) =0.
\eas
The fact that there are no zero divisors in $\bC_q[SU(2)]$ clearly implies that the map is non-degenerate. Note that since $fh^*$ and $h^*f$ are both elements of degree $0$ \wrt to the $\bZ$-grading on $\bC_q[SU_2]$, the image of $g$ does indeed lie in $\bC_q[S^2]$.

% Invertibility
%%%%%%%%%%%%%%%% 
%\bigskip
%%%%%%%%%%%%%%%%

Let us now consider  the element 
\bas
\frak{g} :=  e^+ \oby e^- + q^2 e^- \oby e^+ \in \Om^1_q(S^2) \oby \Om^1_q(S^2).
\eas
It was first considered by Majid in \cite{Maj}, as a $q$-deformation of the standard metric on the two sphere. It has the important property that $\wed(\frak{g}) = 0$, which can be considered as a $q$-deformation of the symmetry of the metric. As is easy to see, $\frak{g}$ is the unique element of  $\Om^1_q(S^2) \oby \Om^1_q(S^2)$ that satisfies the identity 
\bas
(\id, g(\cdot,\w))\frak{g} = \w, & & (g(\cdot,\w),\id)\frak{g} = \w, & & (\text{for all } \w \in \Om^1_q(S^2).
\eas
Motivated by the terminology of \cite{Maj}, we call $\frak{g}$ the {\em inverse} of $\frak{g}$. It will prove important for the definition of the fundamental form of $\frak{g}$.

%\bigskip

Now we can easily extend $g$ to a map from $\Om_q^{\bullet}(S^2)  \oby_{\bC_q[S^2]} \Om_q^{\bullet}(S^2)$ to $\bC_q[S^2]$: For  any $f,g  \in \bC_q[S^2]$,  we set  
\bas
g(f,g) = f g^*, & & g(f e^+ \wed e^-,g e^+ \wed e^-) = fg^*,
\eas 
and moreover require  $\bC_q[S^2], \Om^1_q(S^2)$, and $\Om^2_q(S^2)$ to be orthogonal \wrt $g$.

A {\em Hodge map } $\ast_H$ associated to the metric $(\la \cdot,\cdot \ra)$, was introduced in \cite{LRZGaugedLaplac}. It is defined to be the unique  map $\ast_H: \Om^{\bullet}[\ccp1] \to \Om^{\bullet}[\ccp1]$, for which 
\bas
\la \w, \nu \ra \t = \w \wed (\ast_H(\w')).  
\eas
%uniqueness and existence established using equivariance 
As is easily seen, an explicit description of $\ast$ is given by
\bas
\ast_H(1) = \t, & &  \ast_H(\t) = 1, & &  \asth(e^+) = i e^+, & & \asth(e^-) = - i e^-.  
\eas
As an elementary calculation will verify, $\asth$ commutes with $\ast$.

\section{The Hodge Decompositions}

We will now turn to the question of how to calculate the various cohomology groups of $\Om^\bullet_q(S^2)$. Classically, this is most easily done using Hodge theory, and we shall follow a similar path here. We begin by considering the question of adjoint operators for $\exd, \del$, and $\adel$:

\begin{lem}
The functional $\int$ is closed.
\end{lem}
\demo
Any $\w \in \Om\hol$ can be rewritten as $f \oby e^+$, for some $f \in \E_{-2}$. We then have that 
$
\exd \w = \exd(f \oby e^+) = f\1 \oby \ol{(f\2)^+} \wed e^+
$
Thus, we see that 
\bas
\int \exd(f \oby e^+) = \hat{\t}(h(f\1).\ol{f^+\2} \wed e^+) = \ol{h(f)^+} \wed e^+ = 0
\eas
The proof for $\w' \in \Om\ahol$ is exactly analogous.
\qed

\begin{lem}
With respect to the inner product $\la \cdot, \cdot \ra$, the operators  $\exd, \del$, and $\adel$ are adjointable, with explicit formulae being given by
\begin{align*}
\exd^* = - * \exd *, & & \del^* = - * \ol{\del} *, & & \ol{\del}^* = - * \del *.
\end{align*}
\end{lem}
\demo
The proofs for the three operators are a direct generalisation of the classical proof. For completeness, we present the case for  $\del$: For $\w \in \Om^{k}_q(S^2)$, $\nu \in \Om^{k+1}_q(S^2)$, with $k = 0,1$, we have that 
\bas
\la  \w, \del^* \nu \ra & = - \int g(\w, \ast_H  \adel \ast_H \nu) \t  =  -  \int  \w \wed ( \ast_H^2 (\adel \ast_H \nu)^*) =  (-1)^{k+1} \int  \w \wed ( \del \ast_H \nu^*).\\\
\eas
Now the fact that $\int  \exd =0$, easily implies that $\int \del = 0$, and so, 
\bas                                 
 \int \del \w \wed  \ast_H \nu^*  =  (-1)^{k+1}  \int  \w \wed \del \ast_H \nu^*.
\eas
This tells us that  
\bas
\la  \w, \del^* \nu \ra &  = \int \del \w \wed \ast_H \nu^*  =  \int g(\del \w, \nu)\t  =  \la \del \w, \nu\ra, 
\eas 
as required.
\qed

We call the operator adjoints of $\exd, \del$, and $\adel$, the {\em codifferential}, {\em holomorphic codifferential}, and {\em anti-holomorphic codifferential}, respectively. Using these operators we can introduce $q$-versions of the classical  Dirac and Laplace operators according to 
\bas
D_{\exd} = \exd + \exd^*, & & D_{\del} = \del + \del^*, & & D_{\ol{\del}} = \ol{\del}+\ol{\del}^*,
\eas
and $ \DEL_{\exd} = D_{\exd}^2,  \,\DEL_{\del} = D_{\del}^2$, and  $\DEL_{\adel} = D_{\adel}^2$. 
\begin{lem}
The Hodge map $\ast_H$ commutes with the Laplacians $\lapd,\lapdel$, and $\lapadel$.
\end{lem}
\demo
The proofs are again direct generalisations of the classical versions. For sake of completeness, we present the case for $\del$:
\bas
\asth \circ \lapdel & =  - \asth \del \asth \adel \asth - \asth^2\adel \asth \del =    - \asth \del \asth \adel \asth - \adel \asth \del \asth^2  \\
                            &  =  (- \asth \del \asth \adel \asth - \adel \asth \del \asth) \asth = \lapdel \asth
\eas
\qed

\bigskip

Let us recall the standard basis for $\bC_q[SU(2)]$
$
\{ a^ib^jc^k \,|\, i,j,k \in \bN_0\}.
$
Using it we can define $\bC_q[SU(2)]_k$, for $k \in \bN_0$, as
\bas
\bC_q[SU(2)]_k := \spn_{\bC}\{a^ib^jc^k \,|\, i+j +k = n\}
\eas
and note that each $\bC_q[SU(2)]_k$ is a sub-coalgebra of  $\bC_q[SU(2)]$. We then denote
\bas
\Om_k : = (\bC_q[SU(2)]_k \oby V^{\bullet})^{U_1}
\eas

\begin{lem}
Each  $\Om^1(S^2)_k$ is a left $\bC_q[SU_2]$-comodule, and we have the decomposition 
\bal
\Om^\bullet_q(S^2) = \bigoplus_{k \in \bN_0} \Om_k.
\eal \label{decomp}
Moreover, both $\exd$ and $\exd^*$ restrict to linear endomorphisms of $\Om_k$, for each $k \in \bN_0$.
\end{lem}
\demo
The fact that each $\bC_q[SU(2)]_k$ is a coalgebra directly implies that $\Om^1(S^2)_k$ is a left $\bC_q[SU_2]$-comodule, as a moments thought will confirm.
Now let $\sum_i f_i \oby v_i$ be an element of $\Om^{\bullet}_q(S^2)$, such that $f_i \in \Om_i$, for all $i$, and 
\bas
\DEL_R(\sum_i f_i \oby v_i) = \sum_i f_i \oby v_i  \oby 1.
\eas
Then the fact that  $\Om^1(S^2)_k$ is a left $\bC_q[SU_2]$-comodule implies that 
\bas
\DEL_R(f_i \oby v_i) =  f_i \oby v_i  \oby 1,  & & \text{ for all $i$}.
\eas
The decomposition in (\ref{decomp}) now follows directly.
\qed

With this lemma in hand, we can now move onto the main result of this section.

\begin{thm}
For $\Om^1_q(\Om^\bullet(S^2))$, it holds that:
\begin{enumerate}

\item The Dirac operators $D_{\exd}, D_{\del}$, and $D_{\ol{\del}}$,  as well as the Laplacians $\DEL_{\exd}, \DEL_{\del},$  and $\DEL_{\ol{\del}}$, are diagonalisable;

\item It holds that
\begin{enumerate}
\item $\ker(\DEL_{\exd}) = \ker(D_{\exd}) = \ker(\exd) \cap \ker(\exd^*),$
\item  $\ker(\DEL_{\del}) = \ker(D_{\del}) = \ker(\del) \cap \ker(\del^*),$
\item  $\ker(\DEL_{\ol{\del}}) = \ker(D_{\ol{\del}}) = \ker(\ol{\del}) \cap \ker(\ol{\del}^*)$;
\end{enumerate}

\item We have the three decompositions
\begin{enumerate}
\item $\Om^{\bullet}_q(S^2)  = {\cal H}_{\exd} \oplus \exd( \Om^{\bullet}_q(S^2)) \oplus \exd^*( \Om^{\bullet}_q(S^2))$, 
\item $\Om^{\bullet}_q(S^2)  = {\cal H}_{\del} \oplus \del( \Om^{\bullet}_q(S^2)) \oplus \del^*( \Om^{\bullet}_q(S^2)) $,
 
\item $\Om^{\bullet}_q(S^2)  = {\cal H}_{\ol{\del}} \oplus \ol{\del}( \Om^{\bullet}_q(S^2)) \oplus \ol{\del}^*( \Om^{\bullet}_q(S^2))$.
\end{enumerate}

\end{enumerate}
\end{thm}
\demo 
Let us denote by $\exd_k, \del_k, \ol{\del}, \exd^*_k, \del_k^*$, and $\ol{\del}^*_k$ the respective restrictions  to $\Om_k$ of $\exd, \del, \ol{\del}, \exd^*, \del^*$,  and $\ol{\del}^*$. It is clear that $ \exd^*_k, \del_k^*$, and $\ol{\del}^*_k$ are the adjoints of $\exd_k, \del_k$, and $\ol{\del}_k$ respectively. Since 
\bas
D_{\exd,k} = \exd_k+\exd^*_k, & & D_{\del,k} = \del_k + {\del}_k^*, & & D_{\ol{\del},k} = \ol{\del}_k + \ol{\del}_k^*,
\eas
are all self adjoint operators on a finite dimensional vector space, they are diagonalisable. This  immediately implies that
that 
\bas
\ker(\DEL_{\exd,k}) = \ker(D_{\exd,k}), & & \ker(\DEL_{\del,k}) = \ker(D_{\del,k}), & & \ker(\DEL_{\ol{\del},k}) = \ker(D_{\ol{\del},k}).
\eas

Now since $\la \exd \w, \exd^* \w \ra = \la \exd^2 \w , \nu\ra = 0$, the spaces $\exd(\Om)$ and $\exd^*(\Om)$ are orthogonal. Hence, if $D_\exd(\w)=0$, then 
$d\w=0=d^*\w$.  It follows that $\ker(D)=\ker(d)\cap\ker(d^*)$. That the other versions of the proposition follow is established similarly, and so, we are finished with parts 1 and 2.

\bigskip

Since $\Om_q(S^2)_k$ is a finite dimensional space, we can choose a subspace $X_k \sseq \Om_q(S^2)_k$ \st 
\bas
\Om^{\bullet}_q(S^2)  = {X}_{k} \oplus \exd( \Om^{\bullet}_q(S^2)) \oplus \exd^*( \Om^{\bullet}_q(S^2))
\eas
is an orthogonal decomposition. Now for any $x \in X_k, \w \in \Om^\bullet_q(S^2)$, it must hold that 
\bas
0 = \la x, \exd \w \ra = \la \exd^* x, \w\ra, & & \text{ and }  & & 0 = \la x, \exd^* \w \ra = \la \exd x, \w \ra.
\eas
Hence, we must have that $x \in \ker(\exd) \cap \ker(\exd^*) = {\cal H}_{\exd,k}$, which tells us that $ X_k \sseq {\cal H}_{\exd,k}$. The opposite inclusion  ${\cal H}_{\exd} \sseq X_k$ would follow from the orthogonality of ${\cal H}_{\exd}$ and $ \exd( \Om^{\bullet}_q(S^2)) \oplus \exd^*( \Om^{\bullet}_q(S^2))$. But this is directly implied by   
\bas
\la \exd \w + \exd^* \nu, \mu \ra = \la \exd \w, \mu \ra + \la \exd^* \nu, \mu \ra = 0, & & (\w,\nu \in \Om^\bullet_q(S^2), \mu \in {\cal H}_{\exd}).
\eas
\qed

We call the three decompositions given above {\em Hodge decomposition}, the {\em holomorphic Hodge decomposition}, and the {\em anti-holomorphic Hodge decomposition} respectively. As an easy consequence of the theorem we have following  important result.

\begin{cor}
It holds that 
\bas
\ker(\exd|_{\Om^k})  = {\cal H}_{\exd}^k \oplus \exd( \Om^{k-1}_q(S^2)) \oplus \exd^*( \Om^{k+1}_q(S^2)),\\
\ker(\del|_{\Om^k})  = {\cal H}_{\del}^k \oplus \exd( \Om^{k-1}_q(S^2)) \oplus \exd^*( \Om^{k+1}_q(S^2)),\\
\ker(\ol{\del}|_{\Om^k})  = {\cal H}_{\ol{\del}}^k \oplus \exd( \Om^{k-1}_q(S^2)) \oplus \exd^*( \Om^{k+1}_q(S^2)).
\eas 
and so, we have the following isomorphisms 
\bas
H^k_{\exd} \simeq {\cal H}^k_{\exd}, & & H^{(a,b)}_{\del} \simeq {\cal H}^{(a,b)}_{\exd}, & & H^{(a,b)}_{\exd} \simeq {\cal H}^{(a,b)}_{\exd}.
\eas
\end{cor}
\demo
First note that, since ${\cal H}_{\exd} = \ker(\exd) \cap \ker(\exd^*)$, we cannot have $\exd \w = 0$, for any $\w \in \exd^*(\Om^\bullet_q(S^2)$. This means that 
\bas
\ker(\exd|_{\Om^k})  = {\cal H}_{\exd}^k \oplus \exd( \Om^{k-1}_q(S^2)),
\eas 
which immediately implies that $H^k_{\exd} \simeq {\cal H}^k_{\exd}$. The corresponding isomorphisms for $\del$ and $\ol{\del}$ are established analogously.

\qed

With this result in hand, we can do some cohomological calculations:

\begin{cor}
It holds that 
\begin{enumerate}
\item $H^0 \simeq H^2$, with each having dimension greater than or equal to  $1$;
\item $H_{\del}^{(1,0)} \simeq H^{(0,1)}_{\del}$, and $H_{\adel}^{(1,0)} \simeq H^{(0,1)}_{\adel}$
\end{enumerate}
\end{cor}
\demo
The fact that $H^0 \simeq H^2$ follows directly from the theorem and the fact that the Laplacians commute with Hodge operator. That $H^0$ has dimension greater than or equal to $1$ follows from the basic identity $\DEL_{\exd}(1)  = 0$.%, and that
%\bas
%\DEL_{\exd}(e^+ \wed e^-) & = (\exd \exd^* + \exd^* \exd)(e^+ \wed e^-) = \exd \exd^* (e^+ \wed e^-) \\
%                                                &  = - \exd \circ \ast \circ \exd (1) = 0. 
%\eas

The two isomorphisms in the second part of the theorem are a consequence of the fact that the Laplacians commutate with the $\ast$-operator.
\qed

\section{The Lefschetz Decomposition of $\cp1$}

A major result in the theory of classical Hermitian manifolds is the Lefschetz decomposition. In this section we will formulate a $q$-deformation of the Lefschetz decomposition of the two-sphere.

\bigskip

%%%% kappa
Define the {\em fundamental form} of the metric $g$, using its inverse $\frak{g}$, to be the $2$-form
\bas
\k : = (\id \oby J)\circ \wed (\frak{g}) = -e^+ \wed e^-,
\eas
where $J$ is the almost-complex structure map introduced in the previous section.
It is important to note that $\k$ is a covariant element of $\Om^1_q(S^2)$.
%%%% L
We then define the {Lefschetz operator}
\begin{align*}
L:  \Om^\bullet(\ccp1) \to \Om^\bullet(\ccp1), & & \w \mto \w \wed  \k.
\end{align*}
This operator is covariant due to the fact that $\k$ is coinvariant. 

\bigskip

Classically, the covariance of the Lefschetz operator immediately implies that it is adjointable \wrt the metric $\la \cdot, \cdot \ra$. The following lemma shows that this  carries over to the quantum setting:

\begin{lem}
The operator $L$ has an adjoint $\Lambda$, which we call the {\em dual Lefschetz operator}. Moreover, in direct generalisation of the well-known classical result,  an explicit description of $\Lambda$ is given by
$
\Lambda = \ast \inv \circ L \circ \ast.
$
\end{lem}
\demo
%The proof is simply a case of verifying that the identities in (\ref{adjointability}) hold for all pairs of elements in the set 
%$
%\{f, fe^+,fe^-,fe^+ \wed e^- \,|\, f \in \bC_q[S^2]\}. 
%$
As is easy to see, the only identities which do not hold trivially are 
\bas
\la L f,k e^+ \wed e^-\ra  = \la(f,L(k e^+ \wed e^-)\ra, & &  \text{ and } & & \la L(f e^+ \wed e^-),k\ra = \la fe^+ \wed e^- ,Lk\ra,
\eas
where $h$ is of course another element of $\bC_q[S^2]$. The first identity follows from 
\bas
\la L f,h e^+ \wed e^-\ra & = \int f h^*  = \int g(f,\Lambda(h e^+ \wed e^-)),                                                                                                                                                                                                                                            
\eas
while the second is established similarly. 
\qed

\bigskip

Finally, we introduce the {\em counting operator} on $\Om^{\bullet}_q(S^2)$, defined by setting
\bas
H: \Om_q^{\bullet}(S^2) \to \Om_q^{\bullet}(S^2), & & v \mto \sum_{k=0}^{2} (k-1)\Pi^k,
\eas
where $\Pi^k$ is the projection onto $\Om_q^k(S^2)$. The following lemma relates the counting operator with $L$ and $\Lambda$ in a direct generalisation of a well known classical result:

\begin{lem}
It holds that 
\bas
[H,L] = 2L, & & [H,\Lambda] = -2 \Lambda, & & [L,\Lambda] = H,
\eas
and hence, that the operators $L,\Lambda,$ and $H$ define a representation of $\frak{sl}_2$.
\end{lem}
\demo
We begin by noting that, for any $\w \in \Om^k_q(S^2)$, we have
\bas
[H,L](\w) = (k+1)(\k \wed \w) - \k \wed ((k-1)\w) = 2\k \wed \w,
\eas
which clearly implies the first identity. The second identity follows analogously from 
\bas
[H,\Lambda](\w) = (k-3) \Lambda (\w) - (k-1) \Lambda(\w) = -2 \Lambda (\w).
\eas

The third identity is most easily established in a case by case manner: For $f \in \bC_q[S^2]$, we have 
\bas
[L,\Lambda] f  = L \circ \Lambda (f)  -  \Lambda \circ L (f) = - f  = H(f).
\eas
For $\nu \in \Om^1_q(S^2)$, we have
\bas
[L,\Lambda]\nu & = L \circ \Lambda(\nu) - \Lambda \circ L (\nu)=  \nu = H(\nu).
\eas 
While for $\nu' \in \Om^2(S^2)$, we have
\bas
[L,H]\nu' = L \circ \Lambda(\nu') - \Lambda \circ L(\nu') = 0 = H(\nu').
\eas
\qed

It is interesting to note that this representation of $\frak{sl}_2$ splits into a direct sum of sub-representations according to
\bas
\Om^\bullet_q(S^2) = (\bC_q[S^2] \oplus \Om^2_q(S^2)) \bigoplus \Om^1_q(S^2).
\eas
This is a $q$-deformation of the classical Lefschetz decomposition of the de Rham complex of the two-sphere.

%%%%%%%%%%%%%%%%%%%%%%%%%%%%%%%%%%%%%%%%%%%%%%%%%%%%%%%%%%%%%%%%%%%%%%%%%

\section{The K\"ahler Identities}

A standard result of K\:ahler geometry is that, up to  scalar multiple, these three Laplacians coincide classically. We will now show that this fact carries over to the noncommutative setting by following the standard classical proof based on the K\"ahler identities.

\begin{prop}
We have the following relations:
\begin{align*}
[L,\del^*] = i\ol{\del},        & & [L,\ol{\del}^*]=-i\del,         & & [L,\del] = 0,          & & [L,\ol{\del}] = 0, \\
[\Lambda, \del] = i\ol{\del}^*, & & [\Lambda,\ol{\del}] = -i\del^*, & & [\Lambda, \del^*] = 0, & & [\Lambda,\ol{\del}^*] = 0,
\end{align*}
\end{prop}
\demo
The relations
\begin{align*}
[L,\del] = [L,\ol{\del}] = [\Lambda, \del^*] =  [\Lambda,\ol{\del}^*] = 0
\end{align*}
are direct consequences of the definition of $L$ and $\Lambda$. The remaining relations are easily verified by direct calculation. We show this for the relation $[L,\del^*] = i\ol{\del}$: First we note that $[L,\del^*]$ has a non-zero action only on $\cp1$ and $\Om^{(1,0)}(\ccp1)$. For $f \in \cp1$, we have
\begin{align*}
[L,\del^*] f = (- L \circ * \ol{\del} * + * \ol{\del} * \circ L)f = * \ol{\del} * \circ L f = * \ol{\del} *(f \kappa) = * \ol{\del} f = i \ol{\del} f.
\end{align*}
While for $f\del h \in \Om^{(1,0)}(\ccp1)$, we have
\begin{align*}
[L,\del^*] (f\del h)  & = (- L \circ * \ol{\del} * + * \ol{\del} * \circ L)f\del h = - L \circ * \ol{\del} * f\del h \\
                      & = - i L \circ * \ol{\del}( f\del h).
\end{align*}
Now $\ol{\del}( f\del h) = k e^+ \wed e^-$, for some $k \in \cp1$, and so,
\begin{align*}
[L,\del^*] (f\del h) = - i L \circ * (k e^+ \wed e^-)  - i L (k) = -ik\kappa = i e^+ \wed e^-.
\end{align*}
On the other hand
\begin{align*}
 i\ol{\del} (f\del h) = i k e^+ \wed e^-,
\end{align*}
which establishes the relation.
\qed

\begin{cor}
Denoting by $( \cdot,\cdot)$ the usual anti-commutator bracket, it holds that 
\bas
(\del,\adel^*) = (\adel,\del^*) = 0, & & (\del,\del^*) = (\adel,\adel^*),
\eas
from which follows the identity $\DEL = 2\DEL_{\del} = 2\DEL_{\ol{\del}}$.
\end{cor}
\demo
First note that
\begin{align*}
-i(\ol{\del}\del^* + \del^*\ol{\del}) = \ol{\del}[\Lambda,\ol{\del}]+[\Lambda,\ol{\del}]\ol{\del} = \ol{\del}\Lambda\ol{\del} - \ol{\del}\Lambda \ol{\del} = 0,
\end{align*}
and similarly $\del\ol{\del}^* + \ol{\del}^*\del = 0$. This gives that
\begin{align*}
\DEL & = (\del +  \ol{\del})(\del ^* +  \ol{\del}^*) + (\del^* +  \ol{\del}^*)(\del +  \ol{\del})\\
     & = (\del\del^*+\del^*\del) + (\ol{\del} \ol{\del}^*+ \ol{\del}^* \ol{\del}) + (\ol{\del}\del^* + \del^* \ol{\del}) + (\del \ol{\del}^* +  \ol{\del}^*\del)\\
     & = \DEL_{\del} + \DEL_{\ol{\del}}.
\end{align*}
It remains to show that $ \DEL_{\del} = \DEL_{\ol{\del}}$, which is an easy consequence of the proposition:
\begin{align*}
-i \DEL_{\del} & = -i(\del\del^* + \del^*\del) = \del[\Lambda,\ol{\del}] + [\Lambda,\ol{\del}]\del = \del \Lambda \ol{\del} - \del \ol{\del} \Lambda + \lambda \ol{\del} \del - \ol{\del}\Lambda \del\\
               & = \del\Lambda \ol{\del} + \ol{\del} \del \Lambda - \ol{\del} \Lambda \del = [\del,\Lambda]\ol{\del}+\ol{\del}[\del,\Lambda] = -i\ol{\del}^*\ol{\del}-i\ol{\del}\ol{\del}^* = -i\DEL_{\ol{\del}}.
\end{align*}
\qed

Taken together with Hodge decomposition, this in turn directly implies the following result. It tells us that, just as in the classical case, the Dolbeault cohomology of $\Om^\bullet_q(S^2)$ is a refinement of its de Rham cohomology. 

\begin{cor}
It holds that 
\bas
H^k = \bigoplus_{a+b = k} H_{\del}^{(a,b)} = \bigoplus_{a+b = k} H_{\ol{\del}}^{(a,b)}.
\eas
\end{cor}

%%%%%%%%%%%%%%%%%%%%%%%%%%%%%%%%%%%%%%%%%%%%%%%%%%%%%%%%%%%%%%%%%%%%%%%%%%%%%%%%%%%%%%%%%%%%%%%%%%%%

\bigskip

Mathematical Institute of Charles University, Sokolovsk\'a 83, 186 75 Praha 8, Czech Republic

{\em e-mail}: \tt{obuachalla@karlin.mff.cuni.cz}

\end{document}